\documentclass[12pt]{amsart}
\usepackage{amsfonts, amssymb,amsmath, amsthm, amsxtra, latexsym, amscd}
\usepackage[latin1]{inputenc} 

\newtheorem{theo+}    {Theorem}      [section]
\newtheorem{prop+}  [theo+]  {Proposition}
\newtheorem{coro+}  [theo+]  {Corollary}
\newtheorem{lemm+}  [theo+]  {Lemma}
\newtheorem{deep+}  [theo+]  {Deep Result}
\newtheorem{fact+}  [theo+]  {Fact}
\theoremstyle{definition}
\newtheorem{exam+}  [theo+]  {Example}
\newtheorem{rema+}  [theo+]  {Remark}
\newtheorem{defi+}  [theo+]  {Definition}
\newtheorem{xca+}[theo+]{Exercise}

\numberwithin{equation}{section}

\def\beqn{\begin{equation}}
\def\eeqn{\end{equation}}

\def\epf{\qed \enddemo}

\def\fa{\mathfrak a}

\def\a{\alpha}

\def\Claminv2{|C(\Lambda)|^{-2}}

\def\Ga{\Gamma}
\def\varepsi{\varepsilon}
\def\lam{\lambda}

\def\blam{\underline{\bold \lambda}}

\def\Ome{\Omega}

\def\Aa2D{A^{\a,2}(D)}
\def\bAa2D{\overline{A^{\a,2}(D)}}
\def\Ab2D{A^{\beta,2}(D)}
\def\bAb2D{\overline{A^{\beta,2}(D)}}

\def\Norm#1_#2{\Vert#1\Vert_{#2}}

\def\2pd#1#2{\frac{\partial^2 #1}{\partial #2^2}}
\def\p11d#1#2#3{\frac{\partial^2 #1}{  \partial #2\partial #3  }}

\def\Claminv2{|C(\Lambda)|^{-2}}

\def\varepsi{\varepsilon}
\def\sig{\sigma}

\def\Ga{\Gamma}
\def\lam{\lambda}

\def\tanh{\operatorname{tanh}}

\def\tanh{\operatorname{tanh}}

\def\Aa2D{A^{\a,2}(D)}
\def\bAa2D{\overline{A^{\a,2}(D)}}
\def\Ab2D{A^{\beta,2}(D)}
\def\bAb2D{\overline{A^{\beta,2}(D)}}

\def\m{\underline{\mathbf m}}
\def\n{\underline{\mathbf n}}

\def\ub1#1{\underline{\mathbf 1^{#1}}}

\def\h-g-o-p{hypergeometric orthogonal polynomial }
\def\h-g-o-ps{hypergeometric orthogonal polynomials }

\def\MK{Macdonald-Koornwinder }

\def\bc{\mathbb C}

\def\nat0{\mathbb Z_{\ge 0}} 

\def\bpf{\begin{proof}}
\def\epf{\end{proof}}
\def\beq{\begin{equation}}
\def\eeq{\end{equation}}

\def\bc{\mathbb C}

\def\draft{\centerline{(Draft {\the \day}/{\the\month} \the \year.)}}

\begin{document}

\def\phila{\phi_{\blam}}
\def\vepm{\varepsilon_{\m, \nu}(\blam)}
\def\vepmp{\varepsilon_{\m^\prime, \nu}(\blam)}
\def\kam{\kappa_{\m, \nu}(\blam)}
\def\nutnu{\pi({\nu})\otimes \overline{\pi(\nu)}}
\def\pia2ta2{\pi({\frac a2})\otimes \overline{\pi(\frac a2)}}
\def\Hnu{\mathcal H_{\nu}}
\def\Fnu{\mathcal F_{\nu}}
\def\Pm{\mathcal P_{\m}}
\def\PV{\mathcal P(V_{\bc})}
\def\Pa{\mathcal P(\mathfrak a)}
\def\Pn{\mathcal P_{\n}}
\def\cL#1nu{\mathcal L_{{#1}, \nu}}
\def\cL#1a2{\mathcal L_{{#1}, \frac a2}}
\def\E#1nu{E_{{#1}, \nu}}
\def\E#1a2{E_{{#1}, \frac a2}}

\def\Dc{D_{\mathbb C}}
\def\Vc{V_{\mathbb C}}
\def\Bnu{B_{\nu}} 
\def\el{e_{\blam}} 
\def\bnul{b_{\nu}(\lam)} 
\def\GGa{\Gamma_{\Ome}} 
\def\fc#1#2{\frac{#1}{#2}} 
\def\SS{\mathcal S}

\title[Spherical transform and Jacobi polynomials]
{Spherical transform and Jacobi  polynomials on root systems of type BC}

\author{Genkai Zhang}
\address{Department of Mathematics, Chalmers University of Technology
and G\"o{}teborg University,
S-412 96 G\"o{}teborg, Sweden}
\email{genkai@math.chalmers.se}
\thanks{Research supported by the Swedish Science Research Council (VR)}

\begin{abstract}
Let $R$ be a root system of type BC in $\mathfrak a=\mathbb R^r$ of
general positive multiplicity.
We introduce certain canonical weight function on $\mathbb R^r$
which in the case of symmetric domains
corresponds to the integral kernel
of the Berezin transform. We compute
its spherical transform and 
 prove certain Bernstein-Sato type formula.
This generalizes earlier work of Unterberger-Upmeier,
van Dijk-Pevsner, Neretin  and the author.
Associated to the weight functions there
are Heckman-Opdam orthogonal polynomials of Jacobi type
on the compact torus, after a change of variables 
they form an orthogonal system on the non-compact
space $\mathfrak a$.
We consider their spherical transform
and prove that they are the Macdonald-Koornwinder
polynomials multiplied by the
spherical transform of the canonical weight function. For rank one case this
was proved earlier by Koornwinder.
\end{abstract}

\maketitle
\baselineskip 1.40pc

\section*{Introduction}

The Gaussian functions
and the Hermite polynomials play an important role
in Fourier transform on Euclidean spaces;
the Hermite polynomials diagonalize the
 harmonic oscillator and 
the corresponding Hermite type function
diagonalize the Fourier transform, which
make the Plancherel theory 
more transparent. The generalization
of the Plancherel theory to any non-compact symmetric symmetric
spaces has been studied intensively and
 there are still no general theory
generalizing corresponding the results of  Hermite polynomials
and their Fourier transform, and above all, no concrete
orthogonal systems of  functions with explicit
spherical transforms are constructed. In the present paper we
will establish such a theory for root system of type BC.

Now associated with any root system in $\mathbb R^r$
 there are two
kinds of remarkable orthogonal polynomials,
namely the Heckman-Opdam orthogonal polynomials
giving the spectral decomposition
of the algebra of the Weyl group invariant polynomials
of  Cherednik operators 
acting on certain $L^2$-space on a compact torus
$\mathbb T$,
 and the Macdonald-Koornwinder
polynomials  orthogonal with respect to
certain weight functions on $\mathbb R^r$
defined as a product of Gamma functions.
 Part of the  product 
is in fact the Harish-Chandra Plancherel measure for the
spectral decomposition of the
 algebra of Cherednik operators  on the non-compact space 
$\mathbb R^r$. In the present paper we will
introduce certain canonical weight function  $f_{-2\nu}(t)$
for root system of type BC,
and we compute its
 spherical transform $\widetilde{f_{-2\nu}}(\lam)$. 
 The product of $\widetilde{f_{-2\nu}}(\lam)^2$  with
 the Harish-Chandra measure $|c(\lam)|^{-2}$ gives precisely
weight functions in the Macdonald-Koornwinder orthogonality
relation;  we prove further, roughly speaking,
that the Jacobi polynomials on the compact space, after some change
of variables, and multiplied by the 
function ${f_{-2\nu}}$ 
give orthogonal system of the $L^2$-space on the non-compact
space, and that their spherical transforms are
of the form $\widetilde{f_{-2\nu}}(\lam) p(\lam)$
where $p(\lam)$ are the \MK polynomials.  Thus
the function ${f_{-2\nu}}$ plays the role of the Gaussian
functions in the Fourier analysis whereas
the Jacobi and \MK polynomials play
the role of the Hermite polynomials
and their Fourier transforms (up to 
a multiple of the Gaussian). In rank one case
this has been proved earlier by Koornwinder
\cite{Koornwinder-lnm}.

We remark that
 only the spherical transform
of Weyl group invariant functions
are studied in this paper, yet
we use the tools of the Cherednik operators 
acting on general functions with no invariance,
which makes the computation much easier. Some of our results
can also be easily be generated to that setup.

Our results are motivated by the
study of  the Berezin transform
and branching rule
of holomorphic representations
on weighted Bergman spaces 
on bounded symmetric domains $G/K$ (\cite{gkz-bere-rbsd} and \cite{gz-br2}).
To illustrate our result and explain
some of the motivations we consider
the simplest case of a symmetric
domain, namely the unit disk $G/K=\{z\in \mathbb C; |z|<1\}$
with $G=SU(1, 1)$ and $K=U(1)$.
The root system is $\{\pm 4\varepsilon\}$ on $\fa=\mathbb R$
in our notation in Section 1. The Berezin transform
in question is a group convolution operator 
in the space $L^2(G/K)=L^2(G/K, \frac{dz\wedge d\bar z}{(1-|z|^2)^2})$
 with kernel $(1-|z|^2)^{\nu}$.
The canonical weight function $f_{-2\nu}(t)$
is just the kernel  $(1-|z|^2)^{\nu}$ with $z=\tanh t$ in term of the 
geodesic coordinate $t$. The spherical transform
$\widetilde f_{-2\nu}(\lambda)$
of $f_{-2\nu}(t)$ gives the spectral symbol
of the Berezin transform, and
it has been computed for general bounded symmetric
domain by Unterberger and Upmeier
\cite{UU} (see also
\cite{Dijk-pevzner}, \cite{Neretin-beta-int} and
\cite{gkz-bere-rbsd}). Now there is
an orthogonal basis of the space
$L^2(G/K)^K$ of radial functions
of the form $(1-|z|^2)^{\nu}P(z)$ where
$P(z)$ are the Jacobi  polynomials, $z=\tanh t\in (-1, -1)$.
The Jacobi polynomials (depending on $\nu$) are  
given by the corresponding spherical
polynomials of the same root system with
different multiplicity on the compact torus
 $\mathbb T=i\mathbb R/i\pi\mathbb Z$,
with the change of variable $z=\sin s$
for $s\in  \mathbb T=i\mathbb R/i\pi\mathbb Z$
(see Section 5).  It has been proved
by Koornwinder 
\cite{Koornwinder-lnm} that the spherical
transform of $(1-|z|^2)^{\nu}P(z)$
is of the form $\widetilde f_{-2\nu}(\lambda) Q(\lambda)$
where $Q(\lambda)$ is the Wilson ${}_4F_3$ hypergeometric
orthogonal polynomials. Thus $Q(\lambda)$
are orthogonal with respect to
$\widetilde f_{-2\nu}(\lambda)^2|c(\lambda)|^{-2}d\lambda$.

The \MK polynomials are the multi-variable
generalization of the Wilson's ${}_4F_3$ hypergeometric
orthogonal polynomials,
and can be viewed as the most general
case of a hierarchy of classical polynomials
 \cite{van-Dijen-tams}. We introduce
 the  canonical weight function $f_{-2\nu}$
on a general root system of type BC,
compute the spherical transform
for  functions of the form $f_{-2\nu}P$
with $P$ being the Jacobi polynomials studied
by Heckman-Opdam.

The paper is organized as follows. 
In Section 1 we recall the Plancherel formula
for Opdam-Cherednik transform. In Section 2
we find certain Bernstein-Sato type
formula for the so-called canonical weight
function $f_{\delta}$. Its spherical transform
is computed in Section 3. Finally we compute
the spherical transform of certain Jacobi-type
polynomials in Sections 4 and 5.

It is our  belief that most results in this paper
can be generalized to  general root systems,
that will provide a theory of Hermite-type functions
on non-compact space and thus bring together the orthogonal
polynomials on compact torus and 
spherical transform on non-compact space.

I would like to thank Professors Jacques Faraut, Toshiyuki
Kobayashi, Erik Opdam, Jesper Stokman and Harald Upmeier for
some  helpful discussions.

\section{Spherical transform and Plancherel formula}

Let $\fa=\mathbb R^r$ be an Euclidean space with inner product
$(\cdot, \cdot)$ and let $R\subset \fa^\ast$ be a root system of type BC.
A prototype
of such root systems is 
the restricted root system of  a bounded symmetric domain.
We use some familiar
notation, in order to be consistent with the
notation in the context of bounded symmetric domain
and  also in the context of Plancherel
formula for general root system  \cite{Opdam-acta}.
We fix an orthonormal basis $\{\xi_j\}_{1}^r$ of $\fa$ 
and  a dual basis
$\{\varepsi_j\}$ of $\fa^\ast$, i.e., $\varepsi_j(\xi_k)=\delta_{jk}$,
  so that the positive roots are
$R_+=\{2\varepsi_j; j=1,\cdots, r\}\cup
\{4\varepsi_j; j=1,\cdots, r\}\cup\{2(\varepsi_j\pm \varepsi_k); 1\le
j < k \le r\}$,
with  respective multiplicities $(k_1, k_2, k_3)$
satisfying 
$$2(k_1, k_2, k_3)=(2b, \iota, a).$$
 We assume that $\iota, a, b >0$.
We order the  roots so that
$\varepsi_1 >\varepsi_2 >\cdots > \varepsi_r >0$. The Weyl group
is then $W=S_r\times \mathbb Z_2^r$ consisting of signed permutation
of $\{\xi_j\}$. 
Let $\rho=\frac 12\sum_{\alpha\in R^+ } k_{\alpha}\alpha$
be the half sum of positive roots. Then
$$\rho=\sum_{j=1}^r\rho_j\varepsi_j=\sum_{j=1}^r(\iota +b +(r-j)a)\varepsi_j.
$$

We recall briefly in this section the Plancherel formula
for the spherical transform associated
to the root system $R$ developed by Heckman and Opdam
(\cite{Heckman-Opdam-1},
\cite{Heckman-Opdam-2} and \cite{Opdam-acta}).

Let $d\mu$ be the measure 
$$
d\mu(t)=d\mu_{k}(t)=\prod_{\a\in R_+}|2\sinh (\frac 12\a(t))|^{2k_{\a}}dt
$$
on $\fa$, and let  $L^2(\fa)=L^2(\fa, d\mu)$ and $L^2(\fa)^W$  be the corresponding $L^2$ space
and its  subspace of $W$-invariant functions.

Let 
\begin{equation*}
\begin{split}
D_j&=\partial_j -a\sum_{i<j}\frac{1}{1-e^{-2(t_i-t_j)}}(1-s_{ij})
+ a\sum_{j<k}\frac{1}{1-e^{-2(t_j-t_k)}}(1-s_{jk})+
\\
&\qquad a\sum_{i\ne j}\frac{1}{1-e^{-2(t_j+t_k)}}(1-\sig_{jk})
+2\iota\frac{1}{1-e^{-4t_j}}(1-\sig_{j})
+2b\frac{1}{1-e^{-2t_j}}(1-\sig_{j})
-\rho(\xi_j)
\end{split}
\end{equation*}
be the Cherednik operators acting on functions $f(t)$
on $\fa^{\mathbb C}$, where we identify a function $f(t)$ on $\fa^{\mathbb C}$  with 
$f(t_1, \cdots, t_r)$, for $t=t_1 \xi_1 +\cdots +t_r \xi_r$.
Here $s_{ij}, \sig_{ij}, \sig_i$ are the elements in the 
Weyl group,  $s_{ij}=(ij)$ being the permutation
of $\xi_i$ and $\xi_j$, $\sig_{ij}$ the signed permutation,
$\sig_{ij}(\xi_i)=-\xi_j, \sig_{ij}(\xi_j)=-\xi_i$,
$\sig_{i}$ the reflection $\sig_i(\xi_i)=-\xi_i$,
and all mapping $\xi_k\to \xi_k$ for $k\ne i, j$.
For later purpose we will  also rewrite it as
\begin{equation*}
\begin{split}
D_j&=\partial_j -a\sum_{i<j}\frac{e^{t_i-t_j}}{e^{t_i-t_j}-e^{-(t_i-t_j)}}(1-s_{ij})+
a\sum_{j<k}\frac{e^{t_j-t_k}}{e^{t_j-t_k}-e^{-(t_j-t_k)}}(1-s_{jk})+
\\
&\qquad 
+a\sum_{i\ne j}\frac{e^{t_j+t_k}}{e^{t_j+t_k}-e^{-(t_j+t_k)}}(1-\sig_{jk})
+2\iota\frac{e^{2t_j}}{ e^{2t_j}-e^{-2t_j}}(1-\sig_{j})
\\
&\qquad +2b\frac{e^{t_j}}{e^{t_j}-e^{-t_j}}(1-\sig_{j})
-\rho_j.
\end{split}
\end{equation*}

The operators $\{D_j\}$ are then commuting, and
the  decomposition of $L^2(\fa)$ 
with respect
to the eigenfunctions of $\{D_j\}$
 is given by the Cherednik-Opdam transform \cite{Opdam-acta}, formulated in terms of decomposing representations
of a Hecke algebra. We will mostly be concerned
with the  decomposition of $L^2(\fa)^W$ under $W$-invariant
polynomials of $D_j, j=1, \cdots, r$, which is
given by spherical transform in terms of 
the Heckman-Opdam theory of hypergeometric
functions (\cite{Heckman-Opdam-1}
and \cite{Heckman-Opdam-2}); we shall however use
the operators $D_j$ (in place of their symmetric
polynomials) to
compute the spherical transform of symmetric
functions.

For $\lam\in (\fa^{\ast})^{\mathbb C}$ let $\phi_{\lam}$
be as in \cite{Opdam-acta}
the spherical function. In particular
\begin{equation}\label{D-on-phi}
p(D_1, \cdots, D_r)
\phi_{\lam}=p(\lam(\xi_1), \cdots, \lam(\xi_r))\phi_{\lam},
\end{equation}
for any $W$-invariant polynomial $p$.
 The spherical transform of a function
$f\in L^2(\fa)^W$ is 
$$
\widetilde f(\lam)
=\int_{\fa}f(t)\phi_{\lam}(t)d\mu(t).
$$
The corresponding Plancherel measure
is given by
$$
d\widetilde \mu(\lam)=\frac{(2\pi)^{-r} c_{0}^2}
{ c(\lam) c(-\lam)}d\lam,
$$
where
$$
c(\lam)=\prod_{j=1}^r\frac{  \Gamma(\lam_j +b)  
 \Gamma(2\lam_j)}
{\Gamma(\lam_j +b +\frac{\iota}2) \Gamma(2\lam_j + 2b) }
\prod_{1\le j<k\le r, \epsilon=\pm}^r
\frac{ \Gamma(\lam_j+\epsilon\lam_k )  }
{\Gamma(\lam_j+\epsilon\lam_k +\frac{a}2)}
$$
and
$$
c_0=c(\rho)=\prod_{j=1}^r\frac{  \Gamma(\rho_j +b +1)  
 \Gamma(2\rho_j+1)}
{\Gamma(\rho_j +b +\frac{\iota}2 +1) \Gamma(2\rho_j + b+1) }
\prod_{1\le j<k\le r, \epsilon=\pm}^r\frac{ \Gamma(\rho_j+\epsilon\rho_k +1 )  }
{\Gamma(\rho_j+\epsilon\rho_k +\frac{a}2 +1)}.
$$
Namely, we have
\begin{equation}
\label{plch}
\int_{\fa}|f(t)|^2d\mu(t)
=\int_{i\fa^\ast}|\widetilde f(\lam)|^2 d\widetilde \mu(\lam).
\end{equation}

\section{Bernstein-Sato type formula
for the function $f_\delta$}

We define the weight function $f_{\delta}$
on $\fa$ by
$$f_{\delta}(t)=\prod_{j=1}^r \cosh^\delta t_j.$$
Motivated by the Berezin transform (\cite{UU},
\cite{Dijk-pevzner},
\cite{gkz-bere-rbsd}) we call $f_\delta$
the canonical function.
In the
case of bounded symmetric domains $f_{\delta}$ is the integral
kernel of the Berezin transform considered as a convolution
operator defining the so-called canonical representations, and
is the   analogue of the Gaussian functions in the Euclidean
space, see loc. cit..

In this section we prove the following 
Bernstein-Sato type formula:

\begin{theo+}\label{BSF} There is a $W$-invariant
polynomials of the Cherednik operators $D_j$ mapping
the canonical function $f_\delta$ to $f_{\delta-2}$; more precisely, we have
\begin{equation}
  \label{eq:thm2.1}
\prod_{j=1}^r\left(D_j^2-(\delta+\rho(\xi_1))^2\right)
f_\delta
=\prod_{j=1}^r \left(\delta +a(j-1) \right)
\left( 1-\delta-\iota-a(r-j) \right)
f_{\delta-2}  
\end{equation}
\end{theo+}

The proof of it will be divided into 
two technical lemmas.

\begin{lemm+}\label{lem1}
 The following
formula holds
$$
\prod_{l=1}^j \left(D_l+\delta+\rho(\xi_1)\right)
f_\delta (t)
=\prod_{l=1}^j(\delta+a(l-1)) 
f_\delta (t)
\prod_{l=1}^j(1+\tanh t_l)
$$
\end{lemm+}

\begin{proof} We prove the lemma by induction on $j$. First
we have, since $f_{\delta}$ is Weyl group invariant,
$$
D_1 f_\delta =\delta \tanh t_1  f_\delta -\rho(\xi_1) 
f_\delta,
$$
where we use $\frac{d}{dt}\cosh t = \tanh t \cosh t$.
Rewriting,
$$
(D_1 +\delta +\rho(\xi_1)) f_\delta =\delta (1+\tanh t_1)  f_\delta,
$$
which is the claim for $j=1$. Assume that the equality is true  for 
$\prod_{l=1}^{j-1} (D_l+(\delta+\rho(\xi_1)))$.
We consider it for $j$ in place of $j-1$.
 We need to compute
the operator $D_{j}+\delta+\rho(\xi_{1})$
on the function $
f_\delta  \prod_{l=1}^{j-1}(1+\tanh t_l)
$.
This function is invariant under the permutations
of the first $j-1$ coordinates and the
signed permutations of the last $r-j$ coordinates. Thus
\begin{equation*}
\begin{split}
&\qquad D_{j}f_\delta\prod_{l\le j-1}(1+\tanh t_l)
\\
&=\delta \tanh t_{j} f_{\delta} \prod_{l\le j-1}(1+\tanh t_l)\\
& \quad -a f_{\delta}\sum_{i <j}
\frac{e^{t_i-t_j}} {e^{t_i-t_j}-e^{-(t_i-t_j)} } 
(\tanh t_i -\tanh t_j) \prod_{l\le j-1, l\ne i}(1+\tanh t_l)
\\
&\quad + a f_{\delta} \sum_{i <j}
\frac{e^{t_i+t_j} } {e^{t_i+t_j}-e^{-(t_i+t_j)} } 
(\tanh t_i +\tanh t_{j})\prod_{l\le j-1, l\ne i}(1+\tanh t_l)\\
&\quad-\rho(\xi_{j}) f_\delta  \prod_{l\le j-1}(1+\tanh t_l).
\end{split}
\end{equation*}
Using the formulas 
$\tanh x \pm\tanh y= \frac{\sinh (x\pm y)}{\cosh x \cosh y}$
and $e^x=(1+\tanh x)\cosh x$ we see that the sum
of the $i$th terms in the two summations is, apart
form the factor $af_{\delta}\prod_{l\le j-1, l\ne i}(1+\tanh t_l)$,
\begin{equation*}
\begin{split}
&\qquad \frac{e^{t_i+t_j} } {e^{t_i+t_j}-e^{-(t_i+t_j)} } 
(\tanh t_i +\tanh t_{j})
-\frac{e^{t_i-t_j} } {e^{t_i-t_j}-e^{-(t_i-t_j)} } 
(\tanh t_i -\tanh t_{j})
\\
&=2(1+\tanh t_i)\tanh t_j,
\end{split}
\end{equation*}
and thus
\begin{equation*}
\begin{split}
&\qquad D_{j}f_{\delta}\prod_{l\le j-1}(1+\tanh t_l) 
\\
&=\left( (\delta +a(j-1))\tanh t_{j} -\rho(\xi_{j})\right)
f_{\delta} \prod_{l\le j-1} (1+\tanh t_{l}).
\end{split}
\end{equation*}
It follows then
that
\begin{equation*}
\begin{split}
&\qquad (D_{j}+\delta+\rho(\xi_{1}))
f_\delta\prod_{l\le j-1}(1+\tanh t_l)\\
&=(\delta + a(j-1))f_\delta\prod_{l\le j}(1+\tanh t_l)
+(\rho(\xi_1)-\rho(\xi_j) -a(j-1))f_\delta\prod_{l\le j-1}(1+\tanh t_l).
\end{split}
\end{equation*}
But $\rho(\xi_j)=\iota +b +a(r-j)$ so the second term
vanishes and this completes the proof.
\end{proof}

\begin{lemm+} 
\label{lem2} 
The following
formula holds
\begin{equation*}
\begin{split}
&\phantom{= \,\,}
\prod_{l\ge j} \left(D_l-(\delta+\rho(\xi_1))\right)
f_{\delta-1} e^{t_1 +\cdots +t_r}
\\
&=\prod_{l\ge j}\left(1-\delta-\iota -a(r-l)\right)
f_{\delta-1} e^{t_1 +\cdots +t_r}
 \prod_{l \ge j}(1-\tanh t_l)
\end{split}
\end{equation*}
\end{lemm+}

Accepting temporarily the Lemma, we prove
 Theorem 
 \ref{BSF}. We write the LHS
as
$$\prod_{l=1}^r (D_l-(\delta+\rho(\xi_1)))
\prod_{l=1}^r (D_l+(\delta+\rho(\xi_1)))f_\delta
$$
Taking   $j=r$ in Lemma \ref{lem1}
 we
see that $\prod_{l=1}^r (D_l+(\delta+\rho(\xi_1)))f_\delta$ is, disregarding the constant,
 given by,
$$
f_\delta
\prod_{l=1}^r(1+\tanh t_j)= f_{\delta-1} e^{t_1 +\cdots +t_r}.
$$
Theorem  \ref{BSF} then follows from
 Lemma \ref{lem2}  
for $j=1$, by using the identity
$$
f_{\delta-1} e^{t_1 +\cdots +t_r}\prod_{j=1}^r(1-\tanh t_j)
=f_{\delta-2}.
$$

It remains to prove the Lemma.
\begin{proof} We first compute  $D_r$ and use induction
backward.
We have
\begin{equation*}
\begin{split}
&\phantom{= \,\,} D_r f_{\delta-1} e^{t_1 +\cdots +t_r}\\
&=f_{\delta-1} \prod_{l\le r}e^{t_l}
+(\delta-1) f_{\delta-1}\tanh t_r   \prod_{l\le r}e^{t_l}\\
&\qquad 
+a f_{\delta-1} \sum_{i <r}\frac{e^{t_i+t_r} } {e^{t_i+t_r}-e^{-(t_i+t_r)} }
(e^{t_i +t_r}- e^{-(t_i +t_r)})\prod_{l\le r, l\ne i, r}e^{t_l}
+ \\&\qquad 
+2\iota f_{\delta-1}
\frac{e^{2t_r} } {e^{2t_r}-e^{-2t_r} } (e^{t_r}-e^{-t_r})
\prod_{l\le r-1}e^{t_l}
+ \\&\qquad  +2b f_{\delta-1}
\frac{e^{t_r} } {e^{t_r}-e^{-t_r} } (e^{t_r}-e^{-t_r})
 \prod_{l\le r-1}e^{t_l}
-\rho(\xi_r)f_{\delta-1} \prod_{l\le r}e^{t_l}.
\end{split}
\end{equation*}
The third term can be simplified as
$$
a(r-1)f_{\delta -1} \prod_{l\le r}e^{t_l},
$$
and the next two terms are 
$$
\iota \tanh t_r f_{\delta -1} \prod_{l\le r}e^{t_l}
+(\iota +2b)f_{\delta -1} \prod_{l\le r}e^{t_l}.
$$
So that the previous formula is then
\begin{equation*}
\begin{split}
&\qquad (\delta-1+\iota) \tanh t_r f_{\delta -1} \prod_{l\le r}e^{t_l}
+(1+\iota +2b +a(r-1)-\rho(\xi_r))f_{\delta -1} \prod_{l\le r}e^{t_l}\\
&=(\delta-1+\iota) \tanh t_r f_{\delta -1} \prod_{l\le r}e^{t_l}
+(1+b+\rho(\xi_1) -\rho(\xi_r))f_{\delta -1} \prod_{l\le r}e^{t_l}\\
\end{split}
\end{equation*}
and consequently
\begin{equation*}
\begin{split}
&\qquad (D_r-(\rho_1 +\delta(\xi_1)))
f_{\delta-1} e^{t_1 +\cdots +t_r}\\
&=(\delta-1+\iota) \tanh t_r f_{\delta -1} \prod_{l\le r}e^{t_l}
+(1+b-\delta-\rho(\xi_r))f_{\delta -1} \prod_{l\le r}e^{t_l}\\
&= (\delta-1+\iota)(\tanh t_r -1) f_{\delta -1} \prod_{l\le r}e^{t_l},
\end{split}
\end{equation*}
since $1+b-\delta-\rho(\xi_r)=-(\delta-1+\iota)$. This is the Lemma
for $j=r$.

Assume the lemma is true
for the action of $\prod_{l\ge j+1} (D_l-(\delta+\rho(\xi_1)))$.
We compute
$$
(D_j-(\delta +\rho(\xi_1))) f_{\delta -1} \prod_{l\ge j+1}(1-\tanh t_l) 
 \prod_{l\le r}e^{t_l}
$$
and find it is
$$
-(1-\tanh t_j)(a(r-j) +\iota +\delta -1)
 f_{\delta -1} \prod_{l\ge j}(1-\tanh t_l) 
 \prod_{l\le r}e^{t_l}
$$
by a straightforward yet tedious computation.
\end{proof}

\begin{rema+}
Consider
the double Hecke algebra generated by the Cherednik
operators, the group ring $\mathbb C[W]$, and $\mathbb C[P]$ of polynomials $e^{\varepsi_j}$ as
multiplication operators on $\mathbb C[P]$.
Some refinements of the above computation
then gives certain commutation relations
of the Cherednik operators with $e^{\varepsi_1+\cdots +\varepsi_r}$,
which might be of independent interests. This Theorem
(or its equivalent version under the spherical transform, see
the proof of Theorem 3.2)
can possibly be  obtained also by using 
intertwining relations of the Cherednik-Opdam transform
with the double Hecke algebras proved by Cherednik \cite{Cherednik-imrn97},
it would involves nevertheless many intriguing computations.
\end{rema+}

\section{Spherical transform of the function $f_\delta$}


In this section we shall use Bernstein-Sato type formula
for the function $f_{-2\nu}$ 
to derive a recursive
formula for its spherical transform
$\widetilde{ f_{-2\nu}}(\lam)$
(for sufficiently
large $\nu$). By using a limit
formula  we then derive
a product and, consequently, a Gamma function
formula for $\widetilde f_{-2\nu}(\lam)$.

Denote
$$
\Gamma_{a}(\sig)=\prod_{j=1}^r
\Gamma(\sig-\frac a2(j-1)),
$$
which in the symmetric domain  case is the Gindikin
Gamma function; it will be used to simplify certain product formulas.

We compute  first some normalization constant, which corresponds to the
spherical transform of $f_{-2\nu}$ at $\lam=\rho$.
\begin{lemm+} Let $\nu >\iota +b +a(r-1)$. 
The integral
$$N_{\nu}=\int_{\fa}f_{-2\nu}(t)d\mu(t)
$$
is given by
\begin{equation*}
  \begin{split}
N_{\nu}&=2^{r(2\iota +2b +a(r-1)}r!
\Gamma_{a}(\frac{\iota+1 +2b +a(r-1)}2)\\
&\qquad \times\frac{\Gamma_{a}(\nu-(\frac a2(r-1) +\iota +b))}
{\Gamma_{a}(\nu+\frac{1-\iota}2)}\prod_{1\le i<j\le r}
\frac{\Ga(\frac {a}2(j-i+1)}
{\Ga(\frac{k_3}2(j-i)}.
  \end{split}
\end{equation*}
\end{lemm+}
\proof  By symmetry we need only to integrate
over all $s=(s_1, \cdots, s_r)$ with
$s_j\ge 0$. Making  the change of 
variables $z_j=\tanh^2t_j$, we find first  that
\begin{equation*}
\begin{split}
N_{\nu}&=2^{(2\iota +2b +a(r-1))r}
\int_{[0, 1]^r}
|\prod_{i<j}(z_i-z_j)^{a}|
\prod_{i=1}^r z_i^{\frac {\iota +2b -1}2}
\prod_{i=1}^r (1-z_i)^{\nu -(1+\iota +b +a(r-1))} dz_1 \cdots dz_r \\
\end{split}
\end{equation*}
which is evaluated by the
beta  integral (see \cite[Ex. 7, Sect. 10, Chapt. VII]{Macd-book})
that 
\begin{equation*}
\begin{split}
&\qquad \int_{[0, 1]^r}
|\prod_{i<j}(t_i-t_j)^{a}|
\prod_{i=1}^r t_i^{\alpha-\frac{a}2(r-1)-1}
\prod_{i=1}^r (1-t_i)^{\beta -\frac{a}2(r-1)-1}
 dt_1 \cdots dt_r\\
&=r!\prod_{1\le i<j\le r}
\frac{\Ga(\frac {a}2(j-i+1)}
{\Ga(\frac{a}2(j-i)}
\frac{\Ga_{k}(\alpha)\Ga_{k}(\beta)}
{\Ga_{k}(\alpha+\beta)}.
\end{split}
\end{equation*}
\endproof

In particular it follows that
$$
\frac{N_{\nu+1}}{N_{\nu }}=\prod_{j=1}^r
\frac{\nu-(\frac a2(r-1) +\iota +b)
-\frac a2(j-1)}
{\nu +\frac{1-\iota}2 -\frac a2(j-1)}.
$$
Changing $j-1$ to $r-j$
and using $\rho(\xi_1)=\iota+b +a(r-1)$,
 we can rewrite it as the
following
\begin{equation}
  \label{eq:recur-N}
\frac{N_{\nu+1}}{N_{\nu}}=  
\prod_{j=1}^r
\frac{\nu-\frac12\rho(\xi_1)-\frac12(\iota +b +a(j-1))}
{\nu +\frac {1-\iota}2 -\frac a2(j-1)}.
\end{equation}

\begin{theo+} Let $\nu >\iota +b +a(r-1)$.  The spherical transform of $f_{-2\nu}
$ is given by
$$
\widetilde{f_{-2\nu}}(\lam)
=N_{\nu}
\prod_{j=1}^r \prod_{\epsilon=\pm} 
\frac{
\Gamma(\nu-\frac 12{\rho(\xi_1)}+\epsilon \frac 12 \lam(\xi_j))}{
\Gamma(\nu-\frac 12{\rho(\xi_1)}+\epsilon \frac 12(\iota+b +a(j-1)))}, \quad \lam \in i \fa^\ast.
$$
\end{theo+}

The following result is elementary and we omit its proof; 
it is proved in  \cite{gkz-manu-mat} in the case
when the root system corresponds to a bounded symmetric
domain.
\begin{lemm+} 
Let $\phi$ be a bounded and continuous function
on $\fa$. Then
$$
\lim_{\nu\to \infty}\frac{1}
{N_{\nu}}
\int_{\fa}f_{-2\nu}(t)\phi(t)d\mu(t)
=\phi(0)
$$
\end{lemm+}

We prove now Theorem 3.2.
\begin{proof} We write  $\beta_{\nu}(\lam)=\frac{\widetilde{f_{-2\nu}}(\lam)}{N_{\nu}}
$ and compute it in terms of $\beta_{\nu+1}(\lam)$.
We note first that it is well-defined for $\lam\in i\mathfrak a^\ast$.
Indeed the spherical function
$\phi_{\lam}(t)$ is a bounded function \cite{Opdam-acta}, and
 the function $f_{-2\nu}$
is in $L^1(\fa^\ast, d\mu)$  by Lemma 3.1.
 We perform
the spherical transform on the identity (\ref{eq:thm2.1}) 
with  $\delta =-2\nu$.
Using (\ref{D-on-phi}) we have, 
$$
\prod_{j=1}^r\left( \lam(\xi_j)^2- (-2\nu+\rho(\xi_1))^2\right)
\beta_{\nu}(\lam)
=\prod_{j=1}^r\left(-2\nu+a(j-1)\right)\left( 1+2\nu-\iota-a(r-j)\right)
\beta_{\nu+1}(\lam),
$$
and
\begin{equation*}
\begin{split}
\beta_{\nu}(\lam)
&=\prod_{j=1}^r
\frac{(-2\nu+a(j-1))( 1+2\nu-\iota-a(r-j))}
{(\lam(\xi_j)^2-(-2\nu+\rho(\xi_1))^2)
}
\beta_{\nu+1}(\lam)\\
&=\prod_{j=1}^r\frac{(\nu-\frac{a}2(j-1))( \nu + \frac{1-\iota}2-\frac{a}2(r-j))}
{\left((\nu-\frac 12\rho(\xi_1) +\frac 12\lam(\xi_j)\right)
\left((\nu-\frac 12\rho(\xi_1) +\frac 12\lam(\xi_j)\right)}
\beta_{\nu+1}(\lam)
\end{split}
\end{equation*}
We write further the denominator as
\begin{equation*}
\begin{split}
&\quad\prod_{j=1}^r{(\nu-\frac{a}2(j-1))( \nu + \frac{1-\iota}2-\frac{a}2(r-j))}
\\
&=\prod_{j=1}^r
{(\nu-\frac 12\rho(\xi_1)+\frac12(\iota +b +a(j-1))(\nu + \frac{1-\iota}2-\frac{a}2(j-1))}.
\end{split}
\end{equation*}
Thus
$$
\beta_{\nu}(\lam)
=\frac{{N_{\nu+1}}}{N_{\nu}}
\prod_{j=1}^r
\frac{(\nu-\frac 12\rho(\xi_1)+\frac12(\iota +b +a(j-1)))
(\nu + \frac{1-\iota}2-\frac{a}2(j-1))}
{\left(\nu-\frac 12\rho(\xi_1) +\frac 12\lam(\xi_j)\right)
\left(\nu-\frac 12\rho(\xi_1) -\frac 12\lam(\xi_j)\right)}
\beta_{\nu+1}(\lam).
$$
Using (\ref{eq:recur-N}) this becomes
\begin{equation*}
\begin{split}
\beta_{\nu}(\lam)  
&=\beta_{\nu+1}(\lam)
  \prod_{j=1}^r
\frac{\left(\nu +\frac 12\rho(\xi_1)+\frac12(\iota +b +a(j-1))\right)
\left(\nu +\frac 12\rho(\xi_1)-\frac12(\iota +b +a(j-1))\right)}
{\left(\nu-\frac 12\rho(\xi_1) +\frac 12\lam(\xi_j)\right)
\left(\nu-\frac 12\rho(\xi_1) -\frac 12\lam(\xi_j)\right)}
\\
&=\beta_{\nu+1}(\lam)
\prod_{j=1}^r\left( 1 +\frac{\frac 12(\iota+b +a(j-1))}
{\nu-\frac 12\rho(\xi_1)}\right)
\left( 1 -\frac{\frac 12(\iota+b +a(j-1))}
{\nu-\frac 12\rho(\xi_1)}\right)\\
&\quad \times
\left( 1 + \frac{\frac 12\lam(\xi_j)}{\nu-\frac 12\rho(\xi_1)}\right)^{-1}
\left(1- \frac{\frac 12\lam(\xi_j)}{\nu-\frac 12\rho(\xi_1)}\right)^{-1}.
\end{split}
\end{equation*}
Iterating we find
\begin{equation*}
\begin{split}
\beta_{\nu}(\lam)&=\beta_{\nu+k}(\lam) \prod_{l=1}^{k}\prod_{j=1}^r
\left( 1 +\frac{\frac12(\iota+b +a(j-1))}
{\nu+l-1-\frac 12\rho(\xi_1)}\right)
\left( 1 -\frac{\frac 12(\iota+b +a(j-1))}
{\nu+l-1-\frac 12\rho(\xi_1)}\right)
\\
&\quad \times 
\left( 1 + \frac{\frac 12\lam(\xi_j)}{\nu+l-1-\frac 12\rho(\xi_1)}\right)^{-1}
\left(1- \frac{\frac 12\lam(\xi_j)}{\nu+l-1-\frac 12\rho(\xi_1)}\right)^{-1}
\end{split}
\end{equation*}
However $\frac{1}{N_{\nu+k}} b_{\nu+k}(\lam)\to 1$ as $k\to \infty$
by Lemma 3.3 since the function
$\phi_{\lam}(t)$
is a continuous bounded function
for $\lam\in i\fa^\ast$. We get thus
\begin{equation*}
\begin{split}
\beta_{\nu}(\lam)
&=\prod_{l=1}^{\infty}
\prod_{j=1}^r\left( 1 +\frac{\frac 12(\iota+b +a(j-1))}
{\nu+l-1-\frac 12\rho(\xi_1)}\right)
\left( 1 -\frac{\frac 12(\iota+b +a(j-1))}
{\nu+l-1-\frac 12\rho(\xi_1)}\right)
\\
&\quad\times \left( 1 + \frac{\frac 12\lam(\xi_j)}{\nu+l-1-\frac 12\rho(\xi_1)}\right)^{-1}
\left(1- \frac{\frac 12\lam(\xi_j)}{\nu+l-1-\frac 12\rho(\xi_1)}\right)^{-1}.
\end{split}
\end{equation*}
Our result follows 
rewriting the
infinite product in terms of the Gamma function, 
using
$$
\frac{\Ga(A)\Ga(B)}{\Ga(A+C)\Ga(B-C)}= \prod_{s=0}^{\infty}(1 +\frac C{A+s} )(1-\frac C{B+s});
$$
see
\cite[p.5]{Erdelyi-1}.
\end{proof}

\begin{rema+}
The function $\widetilde {f_{-2\nu}}(\lam)$ for $\nu $ as above
 has exponential decay and  is a bounded function for $\lam\in i\fa^\ast$, and
holomorphic in a strip around $i\fa^\ast$.
We can then define
the Berezin transform on $L^2(\fa, d\mu)^W$
by the inverse spherical transform, namely
$$
B_\nu F(t)=\int_{i\fa^\ast} \widetilde {f_{-2\nu}}(\lam) 
\widetilde F(\lam)\phi_{\lam}(t)
|c(\lam)|^{-2}d\lam.
$$
$B_\nu$ is  then  a bounded and positive  operator
and  it has 
an integral kernel $B(t, s)$. In particular
$B(t, 0)=f_{-2\nu}(t)$. It follows then
from the selfadjoint property of $B$ that
$B(x, y)=f_{-2\nu}(t)f_{-2\nu}(s)L(t, s)$.
It would be interesting
to find a series expansion
for the kernel $L(t, s)$ in terms of the Jack symmetric
polynomials (in the variables $(\tanh^2 t_1, \cdots, \tanh^2 t_r)$).
In the case of bounded symmetric domain, this
is indeed possible, and some degenerate cases of the expansion
have been studied in \cite{gz-br2} (see also
\cite{Beerends-Opdam} for related problems).
\end{rema+}

\section{Spherical transform of a class of functions}

We shall compute the spherical transform
of a class of functions of the form $$
f_{-2\nu}(t)p(x_1^2, \cdots, x_r^2), \qquad
x_j=\tanh t_j, \quad j=1, \cdots, r,
$$
where $p(x_1^2, \cdots, x_r^2)$ is a symmetric
polynomial
in $(x_1, \cdots, x_r)$. We will prove that
they are of the form $\widetilde{f_{-2\nu}}(\lam)
q(\lam)$ where $q$ is a symmetric polynomial
in $(\lam_1^2, \cdots, \lam_r^2)$. Eventually
we will specify $p$ to be Jacobi-type polynomials
in the next section.

We denote $\mathcal D_j$
 the conjugation of $D_j$ by the canonical
function $f_{\delta}$,
$$
\mathcal D_j=\mathcal D_j^{(\delta)}= f_{-\delta} D_j f_{\delta}, \quad 
\quad j=1, \cdots, r.
$$
The operators $\{\mathcal D_j^{(\delta)}\}$ under the change
of variables $x_j=\tanh t_j$ have the form
\begin{equation*}
  \label{eq:calD}
\begin{split}
\mathcal D_j&=\delta x_j +(1-x_j^2)\partial_j
-\frac a2\sum_{i<j}\frac{1+x_i-x_j -x_i x_j}
{x_i-x_j}(1-s_{ij})\\
&+\frac a2\sum_{j<k}\frac{1+x_j-x_k - x_j x_k}
{x_j-x_k}(1-s_{jk})
+\frac a2\sum_{k\ne j}\frac{1+x_j+x_k +x_k x_j}
{x_j+x_k}(1-\sig_{jk})\\
&+\iota \left( 1+\frac 12 (x_j +\frac{1}{x_j})\right)
(1-\sig_j)
+b (1 +\frac{1}{x_j})(1-\sig_j)
-\rho_j  
\end{split}
\end{equation*}

We shall consider the operators $\{\mathcal D_j\}$
acting on polynomials in $x$. For that 
purpose we recall  the  natural ordering $\ge$
on the set of all partitions $\eta=(\eta_1, \cdots, 
\eta_r)\in \mathbb N^r$, $\eta_1\ge \cdots \ge \eta_r\ge 0$. 
For two partitions $\zeta$ and $\eta$,
$\eta\ge \zeta$ if $\eta_1\ge \zeta_1$, $\eta_1+\eta_2\ge \zeta_1+\zeta_2$,
$\cdots$,  $\eta_1+\eta_2+\cdots +\eta_r\ge \zeta_1+\zeta_2+\cdots
+\zeta_r$. For any $\eta\in \mathbb N^r$ denote by $\eta^\ast$the unique partition
in the $S_r$-orbit of $\eta$.

The following lemma can be proved by
a direct elementary computation, it is similar to
the known results that the actions of
the operators $D_j$  on $\mathbb C[P]$
are upper triangular.

\begin{lemm+} The action
of the operators $\mathcal D_j$
 on the monomials $x^\eta =x_1^{\eta_1}\cdots x_r^{\eta_r}$, $\eta\in\mathbb N^r$,
preserves the order in the sense that
$$
\mathcal D_jx^{\eta}= \sum_{\zeta,\, \zeta^\ast \le (\eta^{j})^\ast} a_{\zeta, \eta} x^{\xi},
$$
where
$\eta^{j}= \eta + (0, \cdots, 0, 1, 0, \cdots, 0)$
and the coefficient of $x^{\eta^j}$ is
$$
a_{\eta^{j}, \eta}=\delta + a \#\{i<j: \eta_i >\eta_j\}
+\frac {\iota}2(1-(-1)^{\eta_j}).
$$
\end{lemm+}

\begin{rema+}
Some refinement  of the above lemma can be obtained
by using the Bruhat ordering on the set $\mathbb N^r$, however
we will not need it here.
\end{rema+}

Denoting $m_{\eta}$ the symmetric power sum
\begin{equation}
  \label{eq:sym-sum}
m_{\eta}(y_1, \cdots, y_r)=\sum_{\zeta\in S_{r}\eta} y_1^{\zeta_1}
\cdots y_r^{\zeta_r}.
\end{equation}
We compute now  $m_{\eta}(\mathcal D_1^2, \cdots, \mathcal D_r^2)1$.

\begin{lemm+} The operators
$m_{\eta}(\mathcal D_1^2, \cdots, \mathcal D_r^2) $ acting on the constant
monomial $1$ gives
\begin{equation*}
\quad m_{\eta}(\mathcal D_1^2, \cdots \mathcal D_r^2) 1
=d_{\eta}m_{\eta}(x_1^2, \cdots x_r^2) 
+ \text{lower terms},
\end{equation*}
where the ``lower terms'' here stand for
a linear combination
of symmetric power sums $m_{\eta^\prime}$ of $x_1^2, \cdots, x_r^2$
with $\eta^\prime \prec \eta$,
and
  the leading coefficient is
\begin{equation}
d_{\eta}=\prod_{j=0}^{r}
\prod_{k=0}^{\eta_j-1}(\delta +(r-j)-2k)(\delta +(r-j)a +\iota -1 -2k).
\end{equation}
\end{lemm+}
\begin{proof} This can be
obtained by using the fact that $m_{\eta}(\mathcal D_1^2, \cdots \mathcal D_r^2) 1$
is Weyl group invariant and  by applying Lemma 4.1 successively.
\end{proof}

\begin{prop+} Let $\nu >\iota +b +a(r-1)$ and let $\eta: \eta_1\ge \eta_2
\cdots \ge \eta_r\ge 0$ be a partition.
 The spherical transform
  of the function $f_{-2\nu}(t)m_{\eta}(\tanh^2 t_1, \cdots, \tanh^2 t_r)$ is given by
$$
\widetilde{f_{-2\nu}}(\lam) l_{\eta}(\lam),
$$
and $l_{\eta}(\lam)$ is $W$-invariant polynomial
of $\lam$ with leading term 
 $d_{\eta}^{-1}m_{\eta}(\lam(\xi_1)^2, \cdots,  \lam(\xi_r)^2)$.
\end{prop+}
\begin{proof} We apply Lemma 4.2 for $\delta=-2\nu$. Observe that
$d_{\eta}\ne 0$ for any partition $\eta$. It follows 
that for each partition $\eta$,
$$m_{\eta}(x_1^2, \cdots, x_r^2)
=d_{2\eta}^{-1}m_{\eta}(\mathcal D_1^2, \cdots \mathcal D_r^2) 1
+\sum_{\eta^\prime\prec \eta}
c_{\eta^\prime,\eta}m_{\eta^\prime}(\mathcal D_1^2, \cdots \mathcal D_r^2) 1.
$$
In other words
$$
f_{-2\nu}m_{\eta}(x_1^2, \cdots, x_r^2)
=d_{2\eta}^{-1}m_{\eta}(D_1^2, \cdots, D_r^2) 
f_{-2\nu}
+\sum_{\eta^\prime\prec \eta}
c_{\eta^\prime,\eta}m_{\eta^\prime}(D_1^2, \cdots, D_r^2) f_{-2\nu}.
$$
For each $\zeta$, the spherical transform of 
$m_{\zeta}(D_1^2, \cdots, D_r^2) f_{-2\nu}$
is given by, in view of (\ref{D-on-phi}),
\begin{equation*}
m_{\zeta}(\lam(\xi_1)^2, \cdots, \lam(\xi_r)^2) \widetilde{f_{-2\nu}} (\lam).
\end{equation*}
Our result then follows from the previous formula.
\end{proof}

\section{Jacobi type functions
and \MK polynomials}

We recall briefly the Heckman-Opdam theory
on Jacobi polynomials. We will then construct
certain functions of Jacobi type
and prove that their 
spherical transforms  are Macdonald-Koornwinder type polynomials.

Let $\nu >\iota +b +a(r-1)$.  Consider now the same
root system $R$ with 
$R_+=\{2\varepsi_j; j=1,\cdots, r\}\cup
\{4\varepsi_j; j=1,\cdots, r\}\cup\{2\varepsi_j\pm \varepsi_k; 1\le j< k\le r\}$ 
 and with respective multiplicities
$k^{(\nu)}=(k^{(\nu)}_1, k^{(\nu)}_2, k^{(\nu)}_3)$,
such that
$$
2k^{(\nu)}_2=2(2\nu-(1+\iota +b +a(r-1)))+1, \quad
2(k^{(\nu)}_1+ k^{(\nu)}_2)=\iota +2b, \quad
2k^{(\nu)}_3=a.
$$
(Note that one of the multiplicities, $k^{(\nu)}_1$, is negative,
however they still satisfy the condition $k_{\alpha}+ k_{\frac {\alpha}2}\ge 0$, $\forall \alpha\in R$,
see \cite{Opdam-acta}.)

The corresponding  dual coroot lattice 
is $\frac 12(\mathbb Z\varepsi_1+\cdots 
+\mathbb Z\varepsi_r)=(\frac 12 \mathbb Z)^r$ and
weight lattice is $P=2(\mathbb Z\varepsi_1+\cdots 
+\mathbb Z\varepsi_r)=(2 \mathbb Z)^r$. The set of
dominant weights is 
$$P^+=\{
2(\eta_1\varepsi_1+\cdots 
+\eta_r\varepsi_r); \eta=(\eta_1, \cdots, \eta_r)\in \mathbb N^r, 
\eta_1\ge \cdots \ge \eta_r\ge 0\},$$
and for simplicity we will identify  the elements $2\eta$ in $P^+$
with the partitions $\eta\in \mathbb N^r$, and all $\eta$
appeared below will be assumed to be  partitions.
The partial ordering on the weight lattice $P$
defined in terms of positive roots
is exactly the natural ordering  $\ge$.

The polynomial algebra $\mathbb C[ P]$ of integral
weights is then the polynomial algebra generated
by  $e^{\pm 2\varepsi_j}, j=1, \cdots, r$ as functions
defined on the compact torus  $\mathbb T^r=i\fa/(\pi i\mathbb Z)^r$.
Consider the  inner product in
the space $\mathbb C[ P]$,
\begin{equation}
  \label{eq:inn-on-T}
(f, g)_{k^{(\nu)}}=\int_{[0, \pi]^r} f(s) \overline{g(s)}
|\delta_{k^{(\nu)}}(s)|ds,
\end{equation}
where $\delta_{k^{(\nu)}}(s)=\prod_{\a\in R_+} |2\sin \frac 12\alpha(s))|^{2k^{(\nu)}_{\a}}$.

We recall first the Heckman-Opdam theory
of Jacobi polynomials, see
 \cite{Heckman-Opdam-1}, \cite{Heckman-Opdam-2},
\cite{Opdam-acta}.

\begin{lemm+} For each $\eta=(\eta_1, \cdots, \eta_r)\in \mathbb N^r$, 
$\eta_1\ge \cdots \ge \eta_r\ge 0$,  there exists a unique polynomial
 $P_{\nu, \eta}$
on 
$T$  such that
$$P_{\nu, \eta}=p_{\eta}^W+\sum_{\eta^\prime <\eta} c_{\eta^\prime, \eta}
p_{\eta}^W 
$$
and
$$(P_{\nu, \eta}, p_{\eta^\prime}^W)_{k^{(\nu)}}=0,$$
where $p_{\eta}^W=\sum_{w\in W}w(e^{2(\eta_1\varepsi_1 +\cdots +
\eta_r\varepsi_r)})$ is the Weyl group orbit sum of the power function.
The polynomials $\{P_{\nu, \eta}\}_{\eta}$  forms
an orthogonal basis for $L^2(\mathbb T, |\delta_{k^{(\nu)}}(s)|ds)^W$.
\end{lemm+}
The inner product $\langle P_{\nu, \eta}, P_{\nu, \eta} \rangle_{k^{(\nu)}}$  is
explicitly computed in  \cite{Heckman-Opdam-1} and \cite{Opdam-acta}. 

The Jacobi polynomials on $T$ are symmetric with respect
to the Weyl group $W$, thus they are  symmetric polynomials of
the functions $\frac{e^{2is_j}+e^{-2is_j}}{2}=\cos 2s_j$,
 which can be further identified as
symmetric
polynomials of $(x_1^2, \cdots, x_r^2)$, 
with 
\begin{equation}\label{eq:chg-var}
x_j:=\sin s_j\in [-1, 1], s=s_1\xi_1 +\cdots +s_r\xi_r\in 
\mathbb T=i\mathbb R^r/(\pi i\mathbb Z)^r.
\end{equation}
With some abuse of notations
we denote it also by $P_{\nu, \eta}(x_1, \dots, x_r)$.
We rewrite the characterization of $P_{\nu, \eta}$
in terms of the variables $(x_1, \dots, x_r)$.

The  Jacobi polynomial $P_{\nu, \eta}(x_1, \dots, x_r)$
is then characterized as the unique
polynomial
 so that
$$
P_{\nu, \eta}(x_1, \dots, x_r)=2^{2\eta_1 +\dots +2\eta_r}
m_{\eta}(x_1^2, \dots, x_r^2)
+\sum_{\zeta < \eta} c_{\eta, \zeta}
m_{\zeta}(x_1^2, \dots, x_r^2)
$$
and
that
$$
(P_{\nu, \eta}, m_{\zeta})_{k^{(\nu)}}=0, \quad 
\zeta < \eta
$$
where $(f, g)_{k^{(\nu)}}$ is the
inner product (\ref{eq:inn-on-T}) after the change of variables
 (\ref{eq:chg-var}),
   \begin{equation}
  \label{eq:inn-on-T-var}     
   \begin{split}
(f, g)_{k^{(\nu)}}&=2^{(2\iota +2b +a(r-1))r}\int_{[-1, 1]^r}f(x)\overline{g(x)} \\
&\qquad \times \prod_{i<j}|x_i^2-x_j^2|^{a}
\prod_{j}|x_i^{2b+\iota}| \prod_{j}(1-x_i^2)^{2\nu-(1+\iota +b +a(r-1)}
dx_1\cdots dx_r.
   \end{split}
   \end{equation}

We now define the Jacobi type function on $\fa$, with the
variable $x=\sin s, \, s\in \mathbb T$ in $P_{\nu, \eta}$ being 
$x=\tanh t, \,  t\in \fa$, 
\begin{equation}
  \label{eq:jacobi-type}
H_{\nu, \eta}(t):=f_{-2\nu}(t)P_{\nu, \eta}(\tanh t_1,
\cdots, \tanh t_r).
\end{equation}

\begin{lemm+} The functions $\{H_{\nu, \eta}(t)\}_{\eta}$ form
an orthogonal basis for the space $L^2(\fa)^W$
and 
$$
\langle H_{\nu, \eta}, H_{\nu, \eta}\rangle_{L^2(\fa)}
=\langle P_{\nu, \eta}, P_{\nu, \eta}\rangle_{k^{(\nu)}}
$$
\end{lemm+}
\begin{proof} Changing of variables $t\in \fa \to x\in (-1, 1)^r$,
$x_j:=\tanh t_j$ in the definition of inner product
in $L^2(\fa)$, we  get
the above formula.
The remaining claim follows from (the symmetric version of)
Weierstrass approximation theorem.
\end{proof}

To state the next result we recall briefly (with
some slight reformulation using
the Harish-Chandra $c$-function and  
the function $\widetilde{f_{-2\nu}}(\lam)$)
the \MK polynomials
of type BC (also called multi-variable Wilson polynomials), 
see \cite{van-Dijen-tams}. There exists a unique system
of polynomials $\{p^{MK}_{\eta}(\lam)\}_\eta$ on
$\mathbb R^r$, identified with $i\fa^\ast$, 
such that they are orthogonal
in the space 
$L^2(i\fa^\ast, \widetilde{f_{-2\nu}}(\lam)^2 |c(\lam)|^{-2}d\lam)^W$ and  such that
$p^{MK}_{\eta}(\lam)$ has leading term
$m_{\eta}(\lam(\xi_1)^2, \cdots,  \lam(\xi_r)^2)$.

\begin{theo+}
 The spherical transforms
  of the Jacobi type function $
H_{\nu, \eta}(t)$ on $\fa$
are given by
$$
\widetilde{f_{-2\nu}}(\lam)
q_{\nu, \eta}(\lam), \quad \lam \in i\fa^\ast
$$
where $q_{\nu, \eta}(\lam)$ is a $W$-invariant polynomial.
The polynomials $\{q_{\nu, \eta}(\lam)\}_{\eta}$ form
an orthogonal basis of the space 
$L^2(i\fa^\ast, \widetilde{f_{-2\nu}}(\lam)^2 |c(\lam)|^{-2}d\lam)^W$,
 their norm in that space is given by
$$
\Vert q_{\nu, \eta}\Vert
=\Vert H_{\nu, \eta}\Vert_{L^2(\fa)}
=\Vert P_{\nu, \eta}\Vert_{k^{(\nu)}},
$$
and they are up to constant multiples
the \MK polynomials $\{p^{MK}(\lam)\}$.
\end{theo+}

\begin{proof} That the spherical transform
is of the given form follows from Proposition 4.4, the orthogonality
relation is by the Plancherel formula (\ref{plch}). Now the
polynomials $q_{\eta^\prime }(\lam)$ are orthogonal,
and have the leading term $m_\eta(\lam(\xi_1)^2,
\cdots, \lam(\xi_r)^2)$, again by Proposition 4.4,
thus they are multiples of the 
 \MK
polynomials $p^{MK}(\lam)$, by the uniqueness of the latter 
\cite{van-Dijen-tams}.
 \end{proof}

\newcommand{\noopsort}[1]{} \newcommand{\printfirst}[2]{#1}
  \newcommand{\singleletter}[1]{#1} \newcommand{\switchargs}[2]{#2#1}
  \def\cprime{$'$} \def\cprime{$'$}
\providecommand{\bysame}{\leavevmode\hbox to3em{\hrulefill}\thinspace}
\providecommand{\MR}{\relax\ifhmode\unskip\space\fi MR }
\providecommand{\MRhref}[2]{%
  \href{http://www.ams.org/mathscinet-getitem?mr=#1}{#2}
}
\providecommand{\href}[2]{#2}

\end{document}